\documentclass[12pt,english]{elsarticle}
\usepackage[T1]{fontenc}
\usepackage[latin9]{inputenc}
\usepackage{geometry}
\usepackage{babel}
\usepackage{amsmath}
\usepackage{amssymb}
\usepackage{graphicx}
\usepackage{esint}
\usepackage[unicode=true]
 {hyperref}

\makeatletter

\providecommand{\tabularnewline}{\\}

\newcommand{\Div} {\mathop{\rm div}}

\usepackage[font=footnotesize,labelfont=bf]{caption}
\usepackage{lineno}

\usepackage{sectsty}
\sectionfont{\large}
\subsectionfont{\normalsize\it\bfseries}
\subsubsectionfont{\normalsize\it}

\AtBeginDocument{
  
}

\makeatother

\begin{document}

\title{\textbf{\large{Addressing Integration Error for Polygonal Finite
Elements Through Polynomial Projections: A Patch Test Connection}}}

\author{Cameron Talischi, Glaucio H. Paulino$^{*}$\medskip{}
}

\address{Department of Civil and Environmental Engineering, University of
Illinois at Urbana-Champaign, USA\smallskip{}
}

\address{$^{*}$\emph{Corresponding author, \href{mailto:paulino@uiuc.edu}{paulino@uiuc.edu}}}

\address{\emph{July 16, 2013}}
\begin{abstract}
Polygonal finite elements generally do not pass the patch test as
a result of quadrature error in the evaluation of weak form integrals.
In this work, we examine the consequences of lack of polynomial consistency
and show that it can lead to a deterioration of convergence of the
finite element solutions. We propose a general remedy, inspired by
techniques in the recent literature of mimetic finite differences,
for restoring consistency and thereby ensuring the satisfaction of
the patch test and recovering optimal rates of convergence. The proposed
approach, based on polynomial projections of the basis functions,
allows for the use of moderate number of integration points and brings
the computational cost of polygonal finite elements closer to that
of the commonly used linear triangles and bilinear quadrilaterals.
Numerical studies of a two-dimensional scalar diffusion problem accompany
the theoretical considerations.

\smallskip{}
\end{abstract}
\maketitle

\section{Introduction}

\thispagestyle{plain}

Polygonal finite elements, whose development dates back to the seminal
work of Wachspress \citep{Wachspress}, have gained in popularity
as evidenced by the growing literature on the topic (see, for example,
\citep{Sukumar:2004p401,Sukumar:2006p219,rjasanow2012higher,Gillette2012,PolyTop}).
Among the attractive features of polygonal elements is the greater
flexibility they offer in mesh generation. For example, recently developed
algorithms utilize Voronoi diagrams to generate polygonal and polyhedral
grids with desired regularity and size distribution for complex geometries
\citep{Talischi-PolyMesher,Ebeida,branets2009challenges}. Owing to
their high degree of isotropy, these Voronoi meshes have been recently
used in dynamic fracture simulations to reduce mesh bias in computed
crack patterns \citep{Bishop:2009p2472,SofieDan}. In these analyses,
cracks propagate along element boundaries and commonly-used simplicial
meshes possess preferential crack path directions \citep{Park4K,Rimoli}.
The availability of polygonal finite elements also simplifies mesh
adaption procedures such as local refinement (through element-splitting)
and coarsening (through aggregation) since hanging nodes are naturally
accommodated \citep{Tabarraei:2005p378,Rashid:2006p397,beirao2008posteriori,bassi2012flexibility}.
In several applications, discretization methods on polygonal and polyhedral
grids exhibit improved stability and accuracy when compared to their
simplicial and cubical counterparts. For example, a low-order finite
element scheme defined for large class of polygonal meshes has been
shown to be stable for incompressible flow problems \citep{BeiraoDaVeiga:2010gf,FluidsIJNMF}.
Similarly, polygonal elements can exclude checkerboard layouts and
other numerical instabilities that plague the finite element solution
of topology optimization problems \citep{Talischi:2009p2243,Talischi:2010p3179}.
In terms of accuracy, mixed polygonal finite elements can be more
effective than some commonly-used elements for analysis of incompressible
media \citep{FluidsIJNMF}. Further development of the field can also
contribute to the advancement of compatible or structure-preserving
methods that require computations on dual grids made up of polygonal
and polyhedral cells \citep{hirani2008numerical,codecasa2010new,gillette2011dual,jerome2012analysis}.

The present work deals with the issue of numerical integration for
polygonal finite elements, necessary for the evaluation of weak form
integrals, and its implications for accuracy of the resulting approximations.
Numerical integration for polygonal elements is different from classical
triangular and quadrilateral finite elements for two reasons. First,
few tailored quadrature schemes are available in the literature owing
to the arbitrary geometry of the element domain (see \citep{Mousavi:2009p2586,Natarajan:2009p2477}
and references therein). In practice, a simple but perhaps sub-optimal
procedure is often adopted wherein each polygon is divided into triangular
subdomains and the usual quadrature rules are used in each subdomain.
Second, all the available quadrature schemes are generally inexact
even on regular $n$-gons due to the non-polynomial nature of the
basis functions. One consequence is that the patch test is not passed
unless, of course, a very high-order quadrature rule is used to lower
the errors to machine precision levels. Such a quadrature scheme may
require hundreds of integration points and thus is not feasible \citep{Sukumar:2006p219}.

We will show that the error in the satisfaction of the patch test,
in so far as it measures a lack of polynomial consistency of the discrete
system, places a limit on the convergence of the finite element solutions.
More specifically, the solution error cannot be made smaller with
mesh refinement beyond a certain level, thus rendering the method
non-convergent. A similar issue also plagues meshless methods as they
feature non-polynomial functions and remedies for revival of polynomial
consistency have been explored for quite some time now \citep{chen2001stabilized,puso2008meshfree,liu2010new}.
In a recent series of studies, Babuska, Banerjee and co-workers \citep{Babuska:2008p250,babuvska2009effect,zhang2012numerical}
have shown that, under a zero-sum condition or satisfaction of a discrete
Green's identity, the order of quadrature rule has to be increased
with refinement in order to retain optimal rates of convergence for
meshless discretizations. 

In this work, we consider an alternative approach that ensures the
satisfaction of the patch test and optimal convergence rates with
a \emph{fixed} but sufficiently rich quadrature rule. In practice,
the number of integration points for such a scheme is on the order
of the number of nodes and therefore the overall computational cost
of the method is on par with the linear triangles and quads. As we
will see, we do \emph{not} need to completely eliminate the integration
error in the evaluation of the bilinear form. Instead, what is needed
is the consistency of the discrete bilinear form when one of its arguments
is a piecewise polynomial field. This is accomplished by splitting
the local (elemental) forms according to a polynomial projection of
its arguments and performing numerical integration only on the ``non-polynomial''
part. The subsequent restoration of polynomial consistency is sufficient
 for the satisfaction of the patch test and recovery of optimal convergence
rates. We remark that the present context is somewhat simpler compared
to meshless methods since the the support of basis functions coincides
with element domains and their behavior on element interfaces is known.
As a result, the proposed remedy is carried out at the element level
and directly extended for higher order discretizations.

The proposed approach borrows heavily from the techniques in the mimetic
finite difference (MFD) literature (e.g., \citep{Brezzi:2009ug,DaVeiga:2011p4288})
and in particular the recently developed variational scheme labeled
the Virtual Element Method or VEM \citep{BeiraoDaVeiga:2013,VEM3D}.
The polynomial projection and the splitting of the bilinear form used
in the present work is in fact at the heart of this method. Since
an explicit form of trial and test functions is not available (or
required) in VEM, the remaining non-polynomial term in the bilinear
form is only estimated. What VEM elucidates is a constructive approach
for satisfaction of the patch test, which is sufficient for guaranteeing
the convergence of conforming Galerkin-type approximations. As we
shall see, the linear polygon with the proposed splitting can be in
fact viewed as a particular realization of a first-order Virtual element.
The same cannot be said for higher order elements, including the quadratic
serendipity elements considered here, as the choice of degrees of
freedom will, in general, be different. Similarly, access to the basis
functions greatly simplifies our treatment of non-constant coefficients.
While a thorough comparison between VEM and comparable%
\footnote{A remarkable feature of VEM and related MFD formulations is the systematic
construction of elements with arbitrary order of polynomial accuracy
and continuity on general shapes. Developing comparable finite elements
would be a formidable task.%
} finite elements in terms of cost and accuracy merits its own study,
we emphasize that the goal of the present work is to reduce the burden
of numerical integration for polygonal and polyhedral finite elements,
which we hope, will also be beneficial for a broader class of problems
(e.g. nonlinear problems such as the Navier-Stokes flow \citep{FluidsIJNMF}),
including those for which either an MFD or VEM formulation presently
does not exist.

The remainder of the paper is organized as follows: the model diffusion
problem and its finite element approximation are discussed in the
next section. We consider the relationship between the quadrature
errors and the patch test in section 3. Next, in section 4, we explore
the consequences of the integration error in the convergence of finite
element approximations and discuss sufficient conditions on the discrete
bilinear form to recover optimal convergence rates. We present the
proposed splitting of the bilinear form as well as its constructions
for linear and quadratic elements in section 5. Finally, in section
6, we will address the case of non-constant diffusion tensor. The
paper is concluded with some remarks in section 7 and supplementary
material on construction of polygonal finite elements and implementation
of the proposed approach in the appendix. 

We briefly and partially introduce the notation adopted in this paper.
We denote by $H^{k}(\Omega)$ the standard Sobolev space consisting
of functions whose $k$th derivative is square-integrable over the
given domain $\Omega$ and write $\left\Vert \cdot\right\Vert _{k,\Omega}$
and $\left|\cdot\right|_{k,\Omega}$ for its norm and semi-norm. We
write $L^{2}(\Omega)=H^{0}(\Omega)$ and denote by $H_{g}^{1}(\Omega)$
functions in $H^{1}(\Omega)$ whose trace on $\partial\Omega$ is
equal to $g$. Thus $H_{0}^{1}(\Omega)$ consists of functions that
vanish on the boundary of $\Omega$. For any subset $E\subseteq\Omega$,
we denote by $\left|E\right|$ its Lebesgue measure. The space of
polynomials of degree $m$ over $E$ is denoted by $\mathcal{P}_{m}(E)$.
Finally, an integral evaluated numerically using a quadrature rule
is indicated by $\fint$, assuming that location of integration points
and weights is clear from the context.

\section{Model problem and finite element approximation}

For the sake of concreteness, we focus on a scalar diffusion problem
in two dimensions and limit the discussion to first and second-order
polygonal finite elements. Let $\Omega\subseteq\mathbb{R}^{2}$ be
a bounded open domain with polygonal boundary and consider the steady
state diffusion problem given by
\begin{eqnarray}
-\Div\left(\mathbb{K}\nabla u\right)=f &  & \mbox{in }\Omega\label{eq:StrongForm}\\
u=g &  & \mbox{on }\partial\Omega\label{eq:StrongFormBC}
\end{eqnarray}
where the source $f\in L^{2}(\Omega)$ and boundary data $g\in H^{1/2}(\partial\Omega)$
are prescribed. For the moment, we assume that $\mathbb{K}$ is a
symmetric, positive definite diffusion tensor that is \emph{constant}
over $\Omega$, and postpone the treatment of the general case of
variable coefficients to section 6. 

The variational form of the system of equations (\ref{eq:StrongForm})-(\ref{eq:StrongFormBC})
consists of finding $u\in H_{g}^{1}(\Omega)$ such that
\begin{equation}
a(u,v)=\ell(v),\qquad\forall v\in H_{0}^{1}(\Omega)\label{eq:WeakForm}
\end{equation}
where the bilinear and linear forms are defined, respectively, by
\begin{equation}
a(u,v)=\int_{\Omega}\nabla u\cdot\mathbb{K}\nabla v\mathrm{d}\boldsymbol{x},\qquad\ell(v)=\int_{\Omega}fv\mathrm{d}\boldsymbol{x}
\end{equation}
Note that the coercivity of the bilinear form follows from positive-definiteness
of $\mathbb{K}$ and the fact that $\left|\cdot\right|_{1,\Omega}$
defines a norm on $H_{0}^{1}(\Omega)$. Together with continuity of
the linear form, a consequence of the regularity assumption on $f$,
one can show that the above system has a unique and stable solution.

\subsection{Finite element spaces}

To define the finite element approximation of (\ref{eq:WeakForm}),
we consider a mesh $\mathcal{T}_{h}$ of $\Omega$ consisting of non-overlapping
convex polygons, with maximum diameter $h$. An $H^{1}$-conforming
finite element space associated with the mesh $\mathcal{T}_{h}$ is
given by
\begin{equation}
\mathcal{V}_{h}=\left\{ v_{h}\in C^{0}(\overline{\Omega}):v_{h}|_{E}\in\mathcal{V}_{m}(E),\forall E\in\mathcal{T}_{h}\right\} \label{eq:V_h}
\end{equation}
where $\mathcal{V}_{m}(E)$, is a finite-dimensional space of functions
over element $E$ such that
\begin{equation}
\mathcal{V}_{m}(E)\supseteq\mathcal{P}_{m}(E),\qquad\forall E\in\mathcal{T}_{h}\label{eq:VEsup PE}
\end{equation}
This means that any polynomial of order $m$ can be represented by
the functions in $\mathcal{V}_{m}(E)$. We will consider linear and
quadratic elements that satisfy (\ref{eq:VEsup PE}) with $m=1$ and
$m=2$, respectively. For an $n$-gon $E$, the space $\mathcal{V}_{1}(E)$
has dimension $n$ with degrees of freedom associated with the vertices
of $E$. Similarly, the space $\mathcal{V}_{2}(E)$ is of dimension
$2n$ with additional degrees of freedom associated with mid-points
of each edge.

For the first-order element, the space $\mathcal{V}_{1}(E)$ can be
defined as the span of a set of so-called generalized barycentric
coordinates associated with $E$. A number of such coordinates are
available in the literature \citep{Sukumar:2006p219}, all of which
by definition, are linearly complete and non-negative. The desirable
Lagrangian (Kronecker-delta) property and linear variation on $\partial E$
follow directly from these two properties \citep{Floater:2006p406}. 

In the numerical studies presented in this work%
\footnote{The main results, however, apply to finite elements derived from other
barycentric coordinates (e.g. Mean Value, Sibson, Laplace, maximum
entropy coordinates). %
}, we will use Wachspress coordinates which, under certain shape-regularity
assumptions, yield optimal interpolation estimates. More specifically,
assuming existence of uniform bounds for the aspect ratio, vertex
count and interior angles, we have:
\begin{equation}
\left\Vert u-\pi_{h}u\right\Vert _{1,\Omega}=\mathcal{O}(h)
\end{equation}
for a sufficiently smooth function $u$ \citep{Gillette2012}. In
the above expression, $\pi_{h}u$ denotes the usual nodal interpolation
of $u$ \citep{Brenner-FE}.

For the second-order element, we will use the construction of the
serendipity-like element presented in \citep{rand2011quadratic}.
The basis functions are obtained from appropriate linear combinations
of pairwise products of generalized barycentric coordinates. The resulting
element is constructed to satisfy (\ref{eq:VEsup PE}), the Kronecker-delta
property, and exhibits quadratic variation on the boundary. If Wachspress
coordinates are use for the construction, under the same shape-regularity
assumptions as before, the estimate
\begin{equation}
\left\Vert u-\pi_{h}u\right\Vert _{1,\Omega}=\mathcal{O}(h^{2})
\end{equation}
holds provided that $u$ is sufficiently smooth. Additional details
on the construction of $\mathcal{V}_{m}(E)$ and the polygonal basis
functions can be found in the appendix. If $E$ is a triangle, the
bases reduce to the usual ones and $\mathcal{V}_{m}(E)=\mathcal{P}_{m}(E)$.

\subsection{Approximate problem}

Let $\mathcal{V}_{h,g}=\mathcal{V}_{h}\cap H_{g}^{1}(\Omega)$ and
$\mathcal{V}_{h,0}=\mathcal{V}_{h}\cap H_{0}^{1}(\Omega)$ be the
discrete test and trial spaces%
\footnote{We are tacitly assuming that $g$ and $\mathcal{T}_{h}$ are defined
such that the boundary data can be represented by the trace of functions
in $\mathcal{V}_{h}$ and so $\mathcal{V}_{h,g}$ is non-trivial.
In general, $g$ must be replaced by its nodal interpolation $g_{h}$
and the test space is set to $\mathcal{V}_{h}\cap H_{g_{h}}^{1}(\Omega)$.
However, we will ignore this approximation.%
}. We consider a finite element approximation of (\ref{eq:WeakForm})
that consists of finding $u_{h}\in\mathcal{V}_{h,g}$ such that 
\begin{equation}
a_{h}(u_{h},v_{h})=\ell(v_{h})\qquad\forall v_{h}\in\mathcal{V}_{h,0}\label{eq:GalerkinApprox}
\end{equation}
Here $a_{h}:\mathcal{V}_{h}\times\mathcal{V}_{h}\rightarrow\mathbb{R}$
is a discrete bilinear form defined in terms of symmetric local bilinear
forms $a_{h}^{E}:\mathcal{V}_{m}(E)\times\mathcal{V}_{m}(E)\rightarrow\mathbb{R}$
as follows
\begin{equation}
a_{h}(u,v)=\sum_{E\in\mathcal{T}_{h}}a_{h}^{E}(u,v)
\end{equation}
These local bilinear forms corresponds to the element stiffness matrices
and the summation is related to the assembly process in practice.
For instance, when a quadrature rule is used to compute the stiffness
matrix, we have 
\begin{equation}
a_{h}^{E}(u,v)=\fint_{E}\nabla u\cdot\mathbb{K}\nabla v\mathrm{d}\boldsymbol{x}\label{eq:a_hQ}
\end{equation}
If the quadrature in (\ref{eq:a_hQ}) is exact, $a_{h}(u,v)=a(u,v)$,
and we recover the classical Galerkin approximation. We will consider
alternative constructions of the local bilinear form $a_{h}^{E}$
in section 5. Note that we are assuming in (\ref{eq:GalerkinApprox})
that the exact linear form $\ell$ is available. In practice, numerical
integration is usually used to compute this integral, which amounts
to replacing $\ell(v)$ by 
\begin{equation}
\ell_{h}(v)=\sum_{E\in\mathcal{T}_{h}}\fint_{E}fv\mathrm{d}\boldsymbol{x}\label{eq:lh_Q}
\end{equation}
However, since the main difficulty with numerical integration lies
in the resulting lack of consistency in the bilinear form, we will
ignore this approximation to keep the theoretical discussion simple.
In the motivating examples presented in the next two sections, $f\equiv0$
and numerical integration (\ref{eq:lh_Q}) is in fact exact. Nevertheless
comments will be made throughout regarding the effect of this approximation.

In general, we expect that $a_{h}$ inherits the continuity and coercivity
properties of $a$. These conditions will be satisfied, for example,
if $a_{h}^{E}$ scales as $a^{E}$, that is,%
\footnote{Here $a^{E}$ denotes the restriction of $a$ to element $E$, i.e.,
$a^{E}(u,v)=\int_{E}\nabla u\cdot\mathbb{K}\nabla v\mathrm{d}\boldsymbol{x}$.%
}
\begin{equation}
c_{1}a^{E}(v,v)\leq a_{h}^{E}(v,v)\leq c_{2}a^{E}(v,v)\qquad\forall v\in\mathcal{V}_{m}(E),\ \forall E\in\mathcal{T}_{h}\label{eq:Stability-Cont}
\end{equation}
for some positive constants $c_{1}$ and $c_{2}$, independent of
$h$ and $E$ \citep{BeiraoDaVeiga:2013}. Together with continuity
of $\ell$, we can show that (\ref{eq:GalerkinApprox}) admits a unique
solution $u_{h}$. Additional consistency requirements on $a_{h}$
are naturally needed to ensure convergence of $u_{h}$ to $u$. The
well-known and celebrated patch test provides a means to assess the
consistency of the approximation.

\section{Quadrature error and the patch test}

The engineering patch test is performed by applying boundary conditions
$g=p|_{\partial\Omega}$, with $p\in\mathcal{P}_{m}(\Omega)$, to
a patch of finite elements. This corresponds to the problem where
the exact solution $u=p$. In this section, we consider the approximate
bilinear form defined by numerical integration, cf. (\ref{eq:a_hQ}).

Since $\mathcal{V}_{h}\supseteq\mathcal{P}_{m}(\Omega)$, then $p\in\mathcal{V}_{h,g}$
and we will have $u_{h}=p$ if the quadrature rule in (\ref{eq:a_hQ})
is exact. It can be readily shown that \emph{the patch test is also
passed if} 
\begin{equation}
a_{h}^{E}(p,v)=a^{E}(p,v),\qquad\forall v\in\mathcal{V}_{m}(E),\ \forall E\in\mathcal{T}_{h}\label{eq:ElementConsistency}
\end{equation}
To see this, note that (\ref{eq:ElementConsistency}) implies 
\begin{equation}
a_{h}(p,v_{h})=\sum_{E\in\mathcal{T}_{h}}a_{h}^{E}(p,v_{h})=\sum_{E\in\mathcal{T}_{h}}a^{E}(p,v_{h})=a(p,v_{h})=\ell(v_{h})
\end{equation}
for each $v_{h}\in\mathcal{V}_{h,0}$, and so $u_{h}=p$ is the unique
solution to the discrete problem. \emph{Essentially, (\ref{eq:ElementConsistency})
is a polynomial consistency condition requiring the local bilinear
forms to be exact when one of the arguments is a polynomial function. }

For $m=1$ and the first order patch test, the above condition can
be further simplified. If $\varphi_{1},\dots,\varphi_{n}$ denote
the basis for $\mathcal{V}_{1}(E)$, for (\ref{eq:ElementConsistency})
to hold for\emph{ }an arbitrary linear function $p$, we must have
\begin{equation}
\fint_{E}\nabla\varphi_{i}\mathrm{d}\boldsymbol{x}=\int_{E}\nabla\varphi_{i}\mathrm{d}\boldsymbol{x},\qquad i=1,\dots,n\label{eq:GradCondition}
\end{equation}
Therefore, a sufficient condition for passing the first-order patch
test is that \emph{the quadrature rule integrates the gradient of
the basis functions exactly}. This fact that has been noted and used
in literature of meshless methods (see, for example, \citep{liu2010new}).
\begin{figure}
\begin{centering}
\includegraphics[scale=0.5]{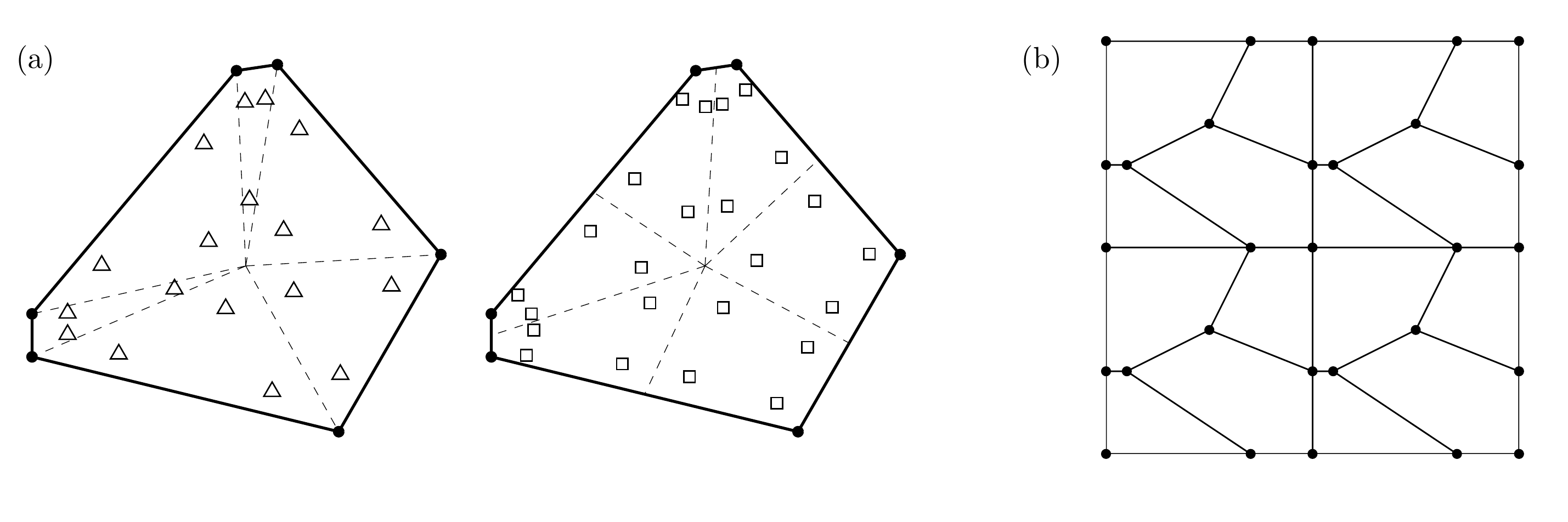}
\par\end{centering}

\caption{(a) Distribution of integration points for second-order triangulation
(left) and quadrangulation (right) integration schemes on a sample
hexagon (b) 2nd-level mesh used for the patch test and convergence
studies. It consists of $2\times2$ patches of two quadrilaterals
and two pentagons.\label{fig:fig0}}
\end{figure}

Few remarks regarding this observation are in order. First, the classical
isoparametric bilinear quadrilateral with $2\times2$ Gauss quadrature
passes the patch test even though there is error in the integration
of discrete bilinear form (i.e., the stiffness matrix) when the elements
are angularly distorted. The patch test is passed precisely because
(\ref{eq:GradCondition}) holds in such a case, a fact seldom discussed
in finite element textbooks. More specifically, if $E$ is the image
of the reference element $\hat{E}=\left[-1,1\right]^{2}$ under the
bilinear map $\boldsymbol{F}$, then the basis functions are defined
through the relation $\varphi_{i}=\hat{\varphi}_{i}\circ\boldsymbol{F}^{-1}$
with $\hat{\varphi}_{i}$ the bilinear function on $\hat{E}$ associated
with $i$th vertex. The relation 
\begin{equation}
\int_{E}\nabla\varphi_{i}\mathrm{d}\boldsymbol{x}=\int_{\hat{E}}\left(D\boldsymbol{F}\right)^{-1}\hat{\nabla}\hat{\varphi}_{i}\det\left(D\boldsymbol{F}\right)\mathrm{d}\hat{\boldsymbol{x}}=\int_{\hat{E}}\mbox{adj}(D\boldsymbol{F})\hat{\nabla}\hat{\varphi}_{i}\mathrm{d}\hat{\boldsymbol{x}}
\end{equation}
indicates the integrand of the right-hand side is bilinear function
of $\hat{\boldsymbol{x}}$ and so $2\times2$ Gauss quadrature on
$\hat{E}$ is exact. In the above expression, $D\boldsymbol{F}$ denotes
the Jacobian matrix for $\boldsymbol{F}$ and $\mbox{adj}(D\boldsymbol{F})$
is the transpose of the cofactor of $D\boldsymbol{F}$. Note, however,
that the bilinear form (\ref{eq:a_hQ}) will not be exact with this
quadrature if $\det\left(D\boldsymbol{F}\right)$ is not constant.

For a general polygon, including distorted quadrilaterals, the Wachspress
basis consisting of rational functions%
\footnote{The basis functions for the iso-parametric quads are also rational
in the physical coordinates but they are images of a polynomial functions
under a polynomial transformation. In fact, the gradient of the Wachspress
basis functions for a general quadrilateral will not be integrated
exactly with Gauss quadrature.%
}, the relation (\ref{eq:GradCondition}) will not hold with the available
quadrature schemes since they are constructed for integration of polynomials.
However, (\ref{eq:GradCondition}) suggests that a quadrature rule
that does a better job of the integrating the gradients of the basis
functions would have a smaller error in the patch test. As mentioned
in the introduction, a simple quadrature scheme for polygonal domains
is obtained by triangulation. We will consider an alternative ``quadrangulation''
procedure, which as shown in Figure \ref{fig:fig0}(a), consists of
splitting the $n$-gon into $n$ quadrilaterals by connecting the
centroid of the polygon to the midpoint of each edge and using Gauss
quadrature (after a bilinear mapping) in each quad. It is evident
from the figure that this approach leads to a denser distribution
of integrations points along the edges of the element where the basis
function gradients are large. By contrast, the triangulation approach
has a denser distribution in the interior of the element.
\begin{table}
\begin{tabular}{|c|c|c|c|c|c|c|}
\hline 
{\small{Quadrature order}} & {\small{1}} & {\small{2}} & {\small{4}} & {\small{8}} & {\small{16}} & {\small{32}}\tabularnewline
\hline 
\hline 
{\small{Triangulation}} & {\small{6.193E-02}} & {\small{1.362E-02}} & {\small{7.165E-04}} & {\small{2.558E-06}} & {\small{4.619E-11}} & {\small{1.223E-15}}\tabularnewline
\hline 
{\small{Quadrangulation}} & {\small{1.766E-02}} & {\small{1.074E-03}} & {\small{5.945E-05}} & {\small{2.206E-07}} & {\small{3.437E-12}} & {\small{1.860E-15}}\tabularnewline
\hline 
\end{tabular}

\caption{Quadrature error for integration of basis function gradients using
different schemes}
\end{table}

In Table 1, we compare the error in the integration of basis function
gradients according to
\begin{equation}
\max_{i}\left|\int_{E}\nabla\varphi_{i}\mathrm{d}\boldsymbol{x}-\fint_{E}\nabla\varphi_{i}\mathrm{d}\boldsymbol{x}\right|
\end{equation}
for the polygon with geometry shown in Figure \ref{fig:fig0}(a).
We first note that the error is finite for all the quadrature orders
considered and many integration points are needed to drive the error
close to machine precision level. Also, the ``quadrangulation''
scheme leads to smaller errors compared to the triangulation approach,
in agreement with the discussion above. We next test to see if the
patch test errors follow the same trend.

\begin{figure}
\begin{centering}
\includegraphics[scale=0.5]{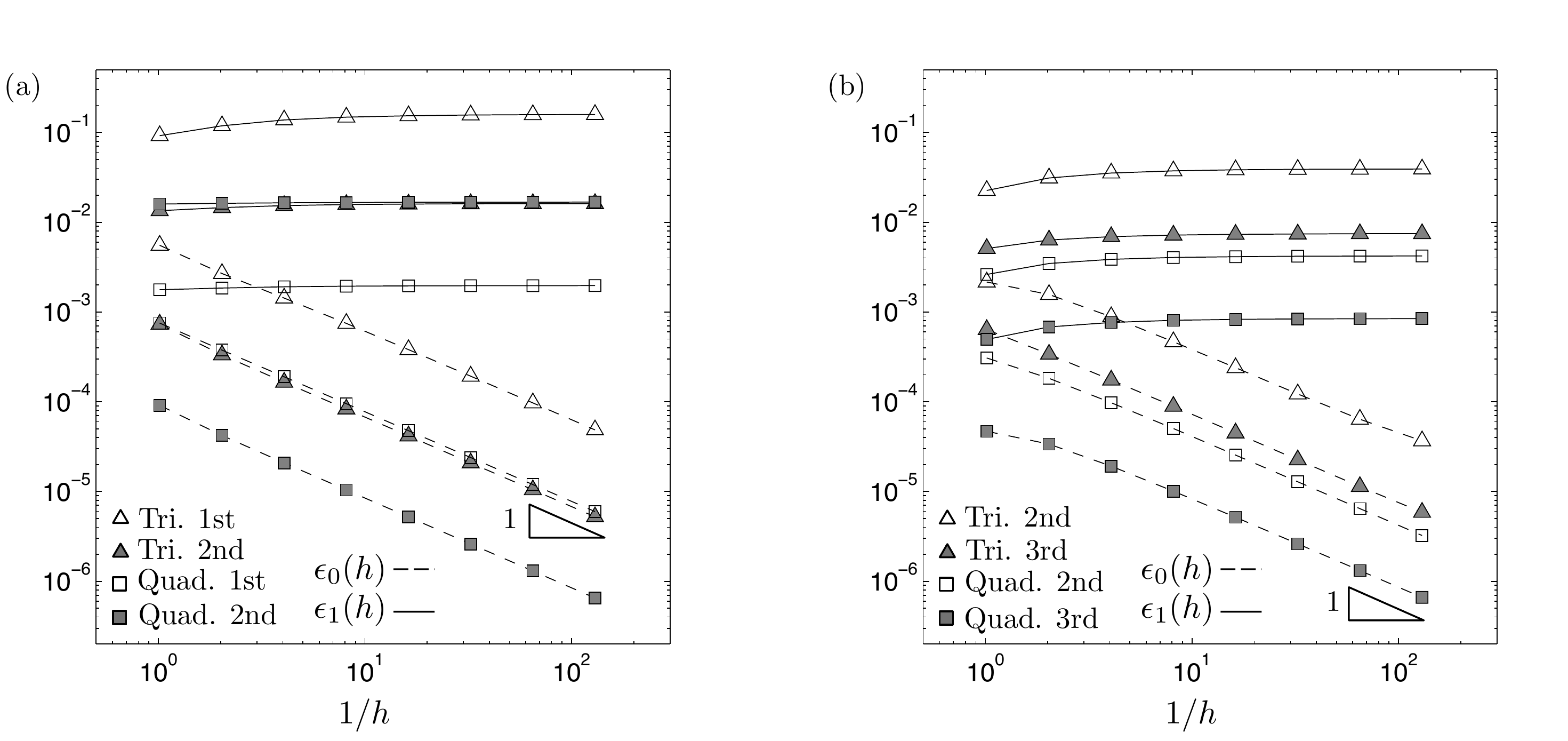}
\par\end{centering}

\centering{}\caption{Results of (a) the linear patch test (b) the quadratic patch test
using indicated quadrature schemes (the legend shows the type of subdivision
and order of quadrature in each subdomain) \label{fig:Fig1}}
\end{figure}

We perform the first-order patch test on the unit square $\Omega=\left]0,1\right[^{2}$
with $u(\boldsymbol{x})=2x_{1}-x_{2}+4$ and, $\mathbb{K}$ taken
to be the identity matrix, on a sequence of polygonal meshes. As shown
in Figure \ref{fig:fig0}(b), the $k$th level mesh consists of $2^{k-1}\times2^{k-1}$
patches of two quadrilaterals and two pentagons. Observe the source
function associated with $u$ vanishes, i.e., $f\equiv-\nabla\cdot(\mathbb{K}\nabla u)=0$.
The reported errors in Figure \ref{fig:Fig1}(a) are with respect
to the $L^{2}$-norm and $H^{1}$-seminorm, given by 
\begin{equation}
\epsilon_{0}(h)=\frac{\left\Vert u-u_{h}\right\Vert _{0,\Omega}}{\left\Vert u\right\Vert _{0,\Omega}},\quad\epsilon_{1}(h)=\frac{\left|u-u_{h}\right|_{1,\Omega}}{\left|u\right|_{1,\Omega}}
\end{equation}
respectively. Note that $\epsilon_{1}(h)$ is the same as the error
in the energy norm for this problem. The results show smaller errors
with the quadrangulation scheme, in agreement with the discussion
above and results of Table 1. We also note\emph{ }that, while the
$L^{2}$-error evidently goes to zero with $h$, the error in the
energy norm does not vanish, indicating that the patch test is not
passed even in a ``weak'' sense. Though not presented, we have observed
the same behavior for higher order quadrature rules.

The above study can be extended to quadratic elements for which condition
(\ref{eq:ElementConsistency}) reduces to a set of conditions on the
integration of basis function gradients, similar to (\ref{eq:GradCondition}),
\emph{and their moments}. We can verify that the errors in the patch
test also correlate with the accuracy of the quadrature scheme for
evaluation of the gradients. 

We will only present the results for a quadratic patch on the same
sequence of meshes as in the previous example. The exact solution
is $u(\boldsymbol{x})=x_{1}^{2}-3x_{1}x_{2}-x_{2}^{2}+5x_{1}$ with
diffusion tensor taken as the identity and $f\equiv0$. Second and
third order quadrature rules are used for the triangular and quadrilateral
subdomains. We observe, from the results shown in Figure \ref{fig:Fig1}(b),
that the quadrangulation scheme again leads to smaller errors that
the triangulation approach. Also, the energy norm errors do not decrease
with mesh refinement, while the $L^{2}$-errors decrease with a linear
rate in the range of mesh sizes considered.

\section{Effects of quadrature error on convergence}

The persistent errors in the patch test under mesh refinement, revealed
in the numerical study of previous section, also indicate the finite
element approximations obtained from (\ref{eq:a_hQ}) using a fixed
quadrature are \emph{not} convergent. Simply put, if solutions do
not converge when the exact solution is a polynomial, the method cannot
be deemed convergent in general. This is an alarming observation and,
to the best of our knowledge, not discussed explicitly in the literature
on polygonal finite elements.

\begin{figure}
\begin{centering}
\includegraphics[scale=0.5]{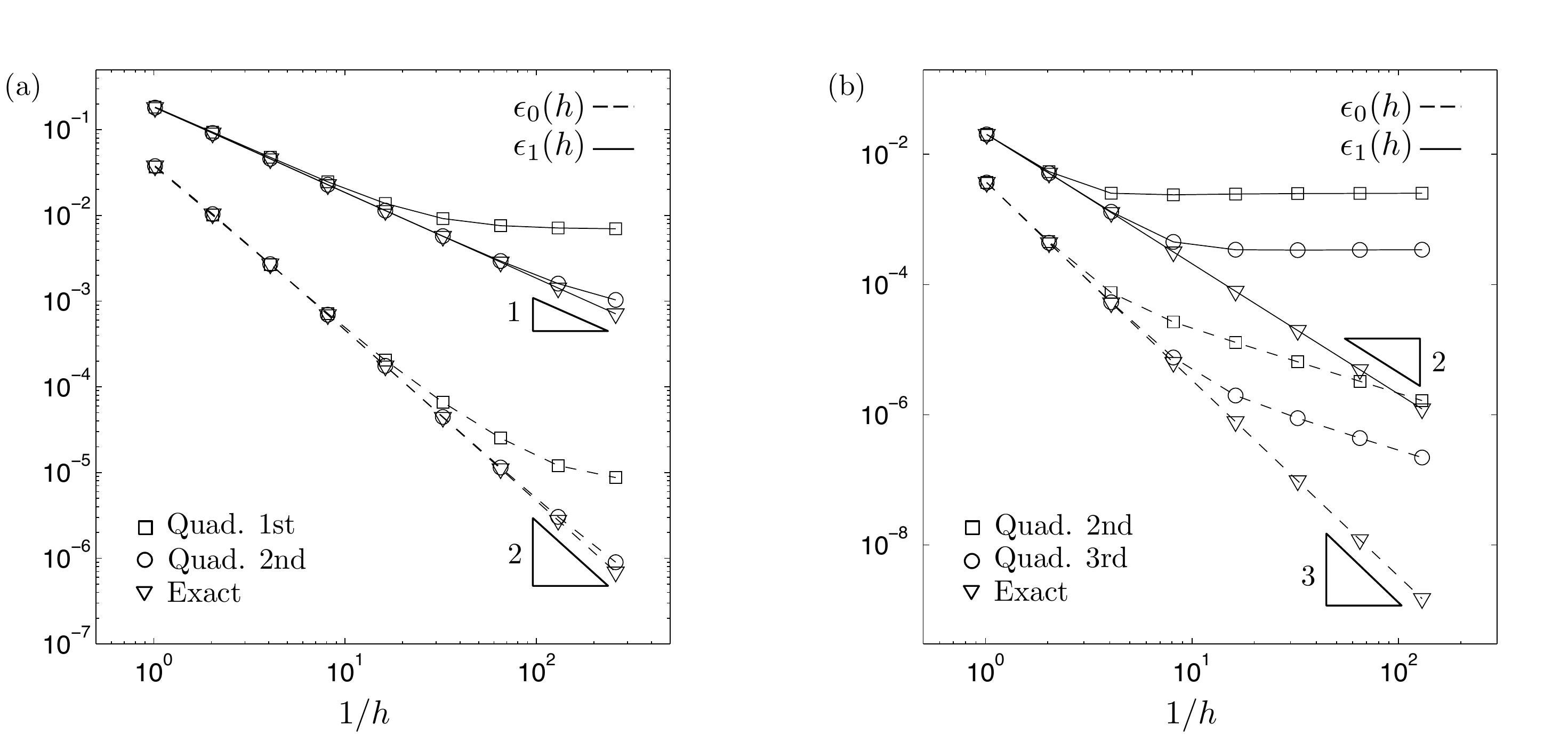}
\par\end{centering}

\caption{Results of the convergence study with (a) linear elements (b) quadratic
elements. In both cases, the quadrangulation scheme with indicated
integration order is used \label{fig:fig2}}
\end{figure}

To further explore the influence of integration error on the convergence
of the approximations, we consider the problem with exact solution
\begin{equation}
u(\boldsymbol{x})=\sin(x_{1})\exp(x_{2})\label{eq:ModelProblem1}
\end{equation}
$\mathbb{K}$ taken as the identity matrix, and $f(\boldsymbol{x})\equiv0$
on the unit square. The boundary data $g$ is specified in accordance
with (\ref{eq:ModelProblem1}). The results for the same regular sequence
of meshes using the quadrangulation scheme is summarized in Figure
\ref{fig:fig2}. While we see optimal rates of convergence with ``exact''
integration%
\footnote{These results are obtained using very high order quadrature such that
error in the calculation of the bilinear form is close to machine
precision levels.%
}, quadrature error leads to a severe degradation of convergence in
the $L^{2}$-norm and a lack of convergence of in the energy norm.
The onset of this poor behavior in fact correlates with the observed
error in the patch test. These results confirm that \emph{the patch
test error places a limit on the accuracy that can be achieved by
the finite element solution. }We remark that this convergence behavior
is reminiscent of error trends for meshless discretizations presented
in reference \citep{zhang2012numerical}. 

We proceed next to discuss a variation of Strang's first lemma (cf.
\citep{Ciarlet:2002:FEM:581834}), that, on the one hand, is in agreement
with the above observation, and, on the other, proves that the satisfaction
of the local consistency condition (\ref{eq:ElementConsistency}),
along with (\ref{eq:Stability-Cont}), is sufficient to ensure optimal
convergence.

As before, let $\pi_{h}$ be the nodal interpolant on $\mathcal{V}_{h}$.
We define $\tau_{h}u$ to be a piecewise polynomial field on $\mathcal{T}_{h}$
that best approximates $u$ with respect to the discrete semi-norm
\begin{equation}
\left|\cdot\right|_{h}^{2}:=\sum_{E\in\mathcal{T}_{h}}\left|\cdot\right|_{1,E}^{2}
\end{equation}
This means that for each $E\in\mathcal{T}_{h}$, the restriction of
$\tau_{h}u$ to $E$ belongs to $\mathcal{P}_{m}(E)$ and
\begin{equation}
\tau_{h}u\bigr|_{E}=\underset{p\in\mathcal{P}_{m}(E)}{\mbox{argmin}}\left|u-p\right|_{1,E}
\end{equation}
Viewed another way, $\left(\nabla\tau_{h}u\right)\bigr|_{E}$ is the
$L^{2}$-projection of $\nabla u$ onto $\left[\mathcal{P}_{m-1}(E)\right]^{2}.$
Note that $\tau_{h}u$ is not necessarily continuous across element
boundaries and need not belong to $H^{1}(\Omega)$. For sufficiently
smooth $u$ and under the shape-regularity assumptions on the elements
in $\mathcal{T}_{h}$, one can show that $\left|u-\tau_{h}u\right|_{h}=\mathcal{O}(h^{m})$
\citep{Brenner-FE}.

Provided that the discrete bilinear form satisfies (\ref{eq:Stability-Cont}),
we have the following \emph{a priori} estimate for the approximation
error:
\begin{equation}
\left|u-u_{h}\right|_{1,\Omega}\leq C\left(\left|u-\pi_{h}u\right|_{1,\Omega}+\left|u-\tau_{h}u\right|_{h}+\sup_{v_{h}\in\mathcal{V}_{h,0}}\sum_{E\in\mathcal{T}_{h}}\frac{\left|a_{h}^{E}(\tau_{h}u,v_{h})-a^{E}(\tau_{h}u,v_{h})\right|}{\left|v_{h}\right|_{1,\Omega}}\right)\label{eq:Strang}
\end{equation}
where $C$ is a constant independent of $h$. The ingredients for
its proof can be found in the proof of Theorem 3.1 in \citep{BeiraoDaVeiga:2013}
and will not be repeated here. In fact, (\ref{eq:Strang}) is at the
core of convergence of VEM.

The first two terms are interpolation errors of order $\mathcal{O}(h^{m})$.
Noting that $\tau_{h}u|_{E}\in\mathcal{P}_{m}(E)$, the remaining
term, involving the approximation of the bilinear form, is closely
related to the errors in respecting the consistency condition (\ref{eq:ElementConsistency}).
If the discrete bilinear form satisfies (\ref{eq:ElementConsistency}),
this term vanishes and one obtains an overall error bound of $\left|u-u_{h}\right|_{1,\Omega}=\mathcal{O}(h^{m})$,
which is optimal. One retains optimal convergence rates even if consistency
error in the approximation of bilinear form is $\mathcal{O}(h^{m})$,
that is, if 
\begin{equation}
\left|a_{h}^{E}(p,v)-a^{E}(p,v)\right|\leq Ch^{m}\left|v\right|_{1,E}\qquad\forall p\in\mathcal{P}_{m}(E),\ \forall v\in\mathcal{V}_{m}(E),\ \forall E\in\mathcal{T}_{h}\label{eq:ElementConsistWeak}
\end{equation}
with the constant $C$ independent of $h$ and $E$.

We note that the estimate (\ref{eq:Strang}) sheds light on convergence
behavior observed in the numerical examples presented in the beginning
of this section. For coarser meshes, the interpolation errors, represented
by the first two terms in (\ref{eq:Strang}), dominate while for sufficiently
small $h$, the consistency error in approximation of the bilinear
form controls the overall error. Thus, the degradation in convergence
is ``delayed'' if the consistency error is lowered. However, the
finite consistency error that accompanies any fixed inexact quadrature
rule will ultimately dominate.

Let us also remark that if the linear form $\ell(v)$ is approximated
through quadrature by (\ref{eq:lh_Q}), an additional term of the
form
\begin{equation}
\sup_{v_{h}\in\mathcal{V}_{h,0}}\frac{\left|\ell_{h}(v_{h})-\ell(v)\right|}{\left|v_{h}\right|_{1,\Omega}}
\end{equation}
will appear in the estimate (\ref{eq:Strang}). However, provided
that the quadrature integrates constant functions exactly on each
element of the mesh and $f$ is sufficiently smooth, the error introduced
is $\mathcal{O}(h^{2})$ (cf. (\ref{eq:QuadraticConsistency}) in
section 5.2) and thus will not affect the rate of convergence of both
linear and quadratic discretizations.

\section{Restoring polynomial consistency}

We now discuss an approach to ensure polynomial consistency even when
using a fixed (but inexact) quadrature rule. The proposed approach
uses a particular representation of the bilinear form $a^{E}$ that
is at central to VEM \citep{BeiraoDaVeiga:2013} and effectively nodal
MFD \citep{Brezzi:2009ug,DaVeiga:2011p4288}. 

Keeping (\ref{eq:VEsup PE}) in mind, we consider a \emph{projection}
operator $\Pi_{m}^{E}:\mathcal{V}_{m}(E)\rightarrow\mathcal{P}_{m}(E)$
for each element $E\in\mathcal{T}_{h}$ such that
\begin{equation}
\begin{cases}
a^{E}(p,\Pi_{m}^{E}v)=a^{E}(p,v)\\
\Pi_{m}^{E}p=p
\end{cases}\quad\forall p\in\mathcal{P}_{m}(E)\label{eq:ProjectionExact}
\end{equation}
Thus $\Pi_{m}^{E}v$ can be thought of as a polynomial approximation
to $v$, as seen by the local bilinear form, minimizing $a^{E}(v-p,v-p)$
in $\mathcal{P}_{m}(E)$. While the above definition is applicable
to other elliptic problems such as elasticity, in the present context
with $\mathbb{K}$ a constant tensor, $\nabla\Pi_{m}^{E}v$ is the
least-squares approximation to $\nabla v$ in $\left[\mathcal{P}_{m-1}(E)\right]^{2}$:
\begin{eqnarray}
\nabla\Pi_{m}^{E}v & = & \underset{\boldsymbol{q}\in\left[\mathcal{P}_{m-1}(E)\right]^{2}}{\mbox{argmin}}\int_{E}\left(\boldsymbol{q}-\nabla v\right)\cdot\mathbb{K}\left(\boldsymbol{q}-\nabla v\right)\mathrm{d}\boldsymbol{x}\nonumber \\
 & = & \underset{\boldsymbol{q}\in\left[\mathcal{P}_{m-1}(E)\right]^{2}}{\mbox{argmin}}\int_{E}\left|\boldsymbol{q}-\nabla v\right|^{2}\mathrm{d}\boldsymbol{x}\label{eq:Pi_LeastSquares}
\end{eqnarray}
Also, we observe that for a triangular element $E$, $\Pi_{m}^{E}$
reduces to the identity map since $\mathcal{V}_{m}(E)=\mathcal{P}_{m}(E)$.

We can use (\ref{eq:ProjectionExact}), along with the symmetry and
linearity of the bilinear form, to show that for $u,v\in\mathcal{V}_{m}(E)$,
that $a^{E}(u,v)$ can be split up as,
\begin{eqnarray}
a^{E}(u,v) & = & a^{E}(\Pi_{m}^{E}u,v)+a^{E}(u-\Pi_{m}^{E}u,v)\nonumber \\
 & = & a^{E}(\Pi_{m}^{E}u,v)+a^{E}(u-\Pi_{m}^{E}u,v)+a^{E}(u-\Pi_{m}^{E}u,\Pi_{m}^{E}v)\nonumber \\
 & = & a^{E}(\Pi_{m}^{E}u,v)+a^{E}(u-\Pi_{m}^{E}u,v-\Pi_{m}^{E}v)\nonumber \\
 & = & a^{E}(\Pi_{m}^{E}u,\Pi_{m}^{E}v)+a^{E}(u-\Pi_{m}^{E}u,v-\Pi_{m}^{E}v)\label{eq:BilinearSplitting}
\end{eqnarray}
Observe that the arguments of the first term are polynomial functions.

Inspired by this identity, we define a discrete bilinear form where
numerical integration is used to evaluate the second ``non-polynomial''
term. That is, we set
\begin{equation}
a_{h}^{E}(u,v):=a^{E}(\Pi_{m}^{E}u,\Pi_{m}^{E}v)+\fint_{E}\nabla\left(u-\Pi_{m}^{E}u\right)\cdot\mathbb{K}\nabla\left(v-\Pi_{m}^{E}v\right)\mathrm{d}\boldsymbol{x}\label{eq:a_hNEW}
\end{equation}
First, we note that once an explicit expression for $\Pi_{m}^{E}$
is derived, the first term can be evaluated exactly because its arguments
are polynomials. Second, this choice of $a_{h}^{E}$ respects the
consistency condition (\ref{eq:ElementConsistency}) since for $u=p\in\mathcal{P}_{m}(E)$,
we have $\Pi_{m}^{E}p=p$ and so 
\begin{equation}
a_{h}^{E}(p,v)=a^{E}(p,\Pi_{m}^{E}v)+\fint_{E}\nabla\left(p-p\right)\cdot\mathbb{K}\nabla\left(v-\Pi_{m}^{E}v\right)\mathrm{d}\boldsymbol{x}=a^{E}(p,\Pi_{m}^{E}v)=a^{E}(p,v)
\end{equation}

The other requirement on the bilinear form, namely condition (\ref{eq:Stability-Cont}),
will be satisfied if a sufficiently rich quadrature scheme is used
for the second term in $a_{h}^{E}$. For example, our numerical studies
confirm that even the lowest order quadrature schemes (triangulation
and quadrangulation) with $n$ integration points are sufficient for
this purpose for the linear elements. In the remainder of this section,
we will discuss how the projection map and the discrete bilinear form
can be computed.

\subsection{Linear elements}

Another key observation made in \citep{BeiraoDaVeiga:2013} is that
the right-hand-side of (\ref{eq:ProjectionExact}) can be computed
exactly given our knowledge of behavior of functions in $\mathcal{V}_{1}(E)$.
Indeed, a simple use of integration by parts shows that for $p\in\mathcal{P}_{1}(E)$
and $v\in\mathcal{V}_{1}(E)$
\begin{eqnarray}
a^{E}(p,v) & = & \int_{E}\nabla v\cdot\mathbb{K}\nabla p\mathrm{d}\boldsymbol{x}\nonumber \\
 & = & -\int_{E}v\Div\left(\mathbb{K}\nabla p\right)\mathrm{d}\boldsymbol{x}+\int_{\partial E}v\mathbb{K}\nabla p\cdot\boldsymbol{n}\mathrm{d}s\label{eq:IntgByParts}\\
 & = & \int_{\partial E}v\mathbb{K}\nabla p\cdot\boldsymbol{n}\mathrm{d}s\nonumber 
\end{eqnarray}
where we have used $\Div\left(\mathbb{K}\nabla p\right)=0$ for the
second equality. Observe that \emph{the last integral can be computed
exactly since $v$ varies linearly on the boundary of $E$.} 

To get an explicit expression for $\Pi_{1}^{E}$, let us set $\boldsymbol{q}=\mathbb{K}\nabla p$
in (\ref{eq:IntgByParts}) to get 
\begin{equation}
\int_{E}\nabla\Pi_{1}^{E}v\cdot\boldsymbol{q}\mathrm{d}\mathbf{x}=a^{E}(p,\Pi_{1}^{E}v)=a^{E}(p,v)=\int_{\partial E}v\boldsymbol{q}\cdot\boldsymbol{n}\mathrm{d}s\label{eq:LinearProj_Step1}
\end{equation}
Because (\ref{eq:ProjectionExact}) must hold for all $p\in\mathcal{P}_{1}(E)$,
we can choose $p$ to recover any arbitrary constant vector $\boldsymbol{q}\in\left[\mathcal{P}_{0}(E)\right]^{2}$
and therefore, (\ref{eq:LinearProj_Step1}) implies, 
\begin{equation}
\int_{E}\nabla\Pi_{1}^{E}v\mathrm{d}\mathbf{x}=\int_{\partial E}v\boldsymbol{n}\mathrm{d}s
\end{equation}
Again observe that the value of the right-hand-side integral depends
only on the\emph{ }nodal values of $v$ and the geometry of $E$.
Moreover, as $\nabla\Pi_{1}^{E}v$ is a constant vector over $E$,
it can be pulled outside of the integral 
\begin{equation}
\nabla\Pi_{1}^{E}v=\frac{1}{\left|E\right|}\int_{\partial E}v\boldsymbol{n}\mathrm{d}s\label{eq:Grad_Pi1}
\end{equation}
This relation could also be seen directly from (\ref{eq:Pi_LeastSquares})
since the best constant approximation to $\nabla v$ over $E$ is
$\left|E\right|^{-1}\int_{E}\nabla v\mathrm{d}\boldsymbol{x}=\left|E\right|^{-1}\int_{\partial E}v\boldsymbol{n}\mathrm{d}s$.

To complete the construction of $\Pi_{1}^{E}$, we assign an appropriate
constant in order to respect the condition $\Pi_{1}^{E}p=p$. We can
choose the constant, for example, for equating $\int_{\partial E}v\mathrm{d}s=\int_{\partial E}\Pi_{1}^{E}v\mathrm{d}s$
or the nodal averages. With the latter choice, we have 
\begin{equation}
\left(\Pi_{1}^{E}v\right)(\boldsymbol{x}):=\overline{v}+\left(\frac{1}{\left|E\right|}\int_{\partial E}v\boldsymbol{n}\mathrm{d}s\right)\cdot(\boldsymbol{x}-\overline{\boldsymbol{x}})\label{eq:LinearProjFinal}
\end{equation}
where the constant $\overline{v}$ is the mean of the nodal values
of $v$ and $\overline{\boldsymbol{x}}$ is the center of $E$ (mean
of the location of vertices of $E$). Clearly the gradient of (\ref{eq:LinearProjFinal})
satisfies (\ref{eq:Grad_Pi1}), and for $p(\boldsymbol{x})=\alpha+\boldsymbol{\beta}\cdot\boldsymbol{x}$,
\begin{eqnarray}
\left(\Pi_{1}^{E}p\right)(\boldsymbol{x}) & = & \overline{p}+\left(\frac{1}{\left|E\right|}\int_{E}\nabla p\mathrm{d}\boldsymbol{x}\right)\cdot\left(\boldsymbol{x}-\overline{\boldsymbol{x}}\right)\nonumber \\
 & = & \left(\alpha+\boldsymbol{\beta}\cdot\overline{\boldsymbol{x}}\right)+\boldsymbol{\beta}\cdot\left(\boldsymbol{x}-\overline{\boldsymbol{x}}\right)\\
 & = & p(\boldsymbol{x})\nonumber 
\end{eqnarray}
verifying that projection map fixes $\mathcal{P}_{1}(E)$.

As a consequence of the form of $\Pi_{1}^{E}$, and the choice of
degrees of freedom for the linear element, the first term in discrete
bilinear form (\ref{eq:a_hNEW}) does not depend on the form of the
basis functions inside the element and is only a function of the geometry
of $E$ and diffusion tensor $\mathbb{K}$. \emph{This means that
elements based on other barycentric coordinates, as well as the first
order VEM formulation, all lead to the same ``consistency'' term.}

We also note that if the quadrature scheme satisfies the gradient
condition (\ref{eq:GradCondition}), then the discrete bilinear form
defined by quadrature (i.e., equation (\ref{eq:a_hQ}) in section
3), is identical to the discrete bilinear form (\ref{eq:a_hNEW}).
Indeed, 
\begin{eqnarray}
a_{h}^{E}(u,v) & = & a^{E}(\Pi_{1}^{E}u,\Pi_{1}^{E}v)+\fint_{E}\nabla\left(u-\Pi_{1}^{E}u\right)\cdot\mathbb{K}\nabla\left(v-\Pi_{1}^{E}v\right)\mathrm{d}\boldsymbol{x}\nonumber \\
 & = & 2a^{E}(\Pi_{1}^{E}u,\Pi_{1}^{E}v)-\fint_{E}\nabla u\cdot\mathbb{K}\nabla\Pi_{1}^{E}v\mathrm{d}\boldsymbol{x}-\fint_{E}\nabla\Pi_{1}^{E}u\cdot\mathbb{K}\nabla v\mathrm{d}\boldsymbol{x}+\fint_{E}\nabla u\cdot\mathbb{K}\nabla v\mathrm{d}\boldsymbol{x}\nonumber \\
 & = & 2a^{E}(\Pi_{1}^{E}u,\Pi_{1}^{E}v)-a^{E}(u,\Pi_{1}^{E}v)-a^{E}(\Pi_{1}^{E}u,v)+\fint_{E}\nabla u\cdot\mathbb{K}\nabla v\mathrm{d}\boldsymbol{x}\nonumber \\
 & = & \fint_{E}\nabla u\cdot\mathbb{K}\nabla v\mathrm{d}\boldsymbol{x}
\end{eqnarray}
This implies that, for isoparameteric bilinear quads, applying quadrature
to either representation of the local bilinear form yields the same
result.

\subsection{Quadratic elements}

For the serendipity element considered here, $a^{E}(p,v)$ with $p\in\mathcal{P}_{2}(E)$
\emph{cannot} be reduced to an integral on the boundary of the element.
Therefore, numerical quadrature will be needed for the evaluation
of an area integral. This can be seen from (\ref{eq:IntgByParts})
and the fact that for $p\in\mathcal{P}_{2}(E)$, the quantity $\Div(\mathbb{K}\nabla p)$
does not necessarily vanish. However, the quadrature error for computing
this term is $\mathcal{O}(h^{2})$ since the integrand contains the
basis functions and not their gradients. This error is sufficient
for ensuring the consistency condition (\ref{eq:ElementConsistWeak})
and subsequently maintaining optimal convergence rates. 

In this case, we define the discrete bilinear to be again given by
(\ref{eq:a_hNEW}), but change the definition of the projection map
to 
\begin{equation}
\begin{cases}
a^{E}(p,\Pi_{2}^{E}v)=-\fint_{E}v\Div\left(\mathbb{K}\nabla p\right)\mathrm{d}\boldsymbol{x}+\int_{\partial E}v\mathbb{K}\nabla p\cdot\boldsymbol{n}\mathrm{d}s\\
\Pi_{2}^{E}p=p
\end{cases}\quad\forall p\in\mathcal{P}_{2}(E)\label{eq:ProjecAppQuad}
\end{equation}
which is a slight deviation from (\ref{eq:ProjectionExact}) with
a revisited right-hand side for the first expression. 

As before, the boundary integral in (\ref{eq:ProjecAppQuad}) can
be computed exactly since the integrand is a polynomial. For the two
conditions to be consistent, that is, for the first equality to hold
when $v=q\in\mathcal{P}_{2}(E)$, \emph{we must require that the quadrature
rule is exact for second order polynomials}, that is,
\begin{equation}
\fint_{E}q\mathrm{d}\boldsymbol{x}=\int_{E}q\mathrm{d}\boldsymbol{x},\qquad\forall q\in\mathcal{P}_{2}(E)
\end{equation}
This, in particular, indicates that the first order triangulation
and quadrangulation schemes consisting of $n$-points will not be
sufficient for the quadratic elements and a second-order accurate
quadrature must be used in the subdomains. These rules have proven
in our numerical studies to also be sufficient for ensuring that condition
(\ref{eq:Stability-Cont}) is met. A possible alternative is to use
the quadrature rules in \citep{Mousavi:2009p2586} which are constructed
for exact integration of polynomials on polygonal domains. Compared
to subdivision schemes considered here, they require fewer evaluation
points to achieve quadratic precision.

To verify the satisfaction of the weaker consistency condition (\ref{eq:ElementConsistWeak}),
let $p\in\mathcal{P}_{2}(E)$ and set $c\equiv\Div\left(\mathbb{K}\nabla p\right)$.
Then, for $v\in\mathcal{V}_{2}(E)$, we have
\begin{eqnarray}
\left|a_{h}^{E}(p,v)-a^{E}(p,v)\right| & = & \left|a^{E}(p,\Pi_{2}^{E}v)-a^{E}(p,v)\right|\nonumber \\
 & = & \left|\int_{E}v\Div\left(\mathbb{K}\nabla p\right)\mathrm{d}\boldsymbol{x}-\fint_{E}v\Div\left(\mathbb{K}\nabla p\right)\mathrm{d}\boldsymbol{x}\right|\label{eq:QuadraticConsistency}\\
 & = & \left|c\right|\left|\int_{E}v\mathrm{d}\boldsymbol{x}-\fint_{E}v\mathrm{d}\boldsymbol{x}\right|\nonumber 
\end{eqnarray}
The last expression is bounded by $Ch^{2}\left|v\right|_{1,2,E}$
provided that the quadrature integrates constant functions exactly
(see Exercise 4.1.4 of \citep{Ciarlet:2002:FEM:581834}). 

The fact that the stronger consistency condition (\ref{eq:ElementConsistency})
is not satisfied for the quadratic element suggests that the patch
test will only be passed asymptotically with mesh refinement. However,
we can directly show that the patch test will be passed \emph{exactly
on any mesh} if the approximate linear form $\ell_{h}$ defined by
(\ref{eq:lh_Q}) is used with the same quadrature rule as that of
the discrete bilinear form. As mentioned before, this is usually the
case in practice. Indeed for $p\in\mathcal{P}_{2}(\Omega)$ and $v_{h}\in\mathcal{V}_{h,0}$,
we have
\begin{eqnarray}
a_{h}(p,v_{h}) & = & \sum_{E\in\mathcal{T}_{h}}a_{h}^{E}(p,v_{h})\nonumber \\
 & = & \sum_{E\in\mathcal{T}_{h}}a^{E}(p,\Pi_{2}^{E}v_{h})\nonumber \\
 & = & \sum_{E\in\mathcal{T}_{h}}-\fint_{E}v_{h}\Div\left(\mathbb{K}\nabla p\right)\mathrm{d}\boldsymbol{x}+\sum_{E\in\mathcal{T}_{h}}\int_{\partial E}v_{h}\mathbb{K}\nabla p\cdot\boldsymbol{n}\mathrm{d}s\label{eq:PatchTestQuad}\\
 & = & \sum_{E\in\mathcal{T}_{h}}\fint_{E}v_{h}f\mathrm{d}\boldsymbol{x}\nonumber \\
 & = & \ell_{h}(v_{h})\nonumber 
\end{eqnarray}
Note that the second term in (\ref{eq:PatchTestQuad}) cancels out
as the internal edges of the mesh are visited twice (the normal vector
$\boldsymbol{n}$ changes sign each time) and $v_{h}=0$ on the boundary
edges. Also, we set $f=-\Div(\mathbb{K}\nabla p)$ in the second to
last equality. Our numerical studies in fact confirm that the quadratic
patch test will be passed up to machine precision errors. 

We have provided the details on an explicit construction of the projection
map for quadratic elements in the appendix. As for the linear element,
it is completely characterized by the two conditions in (\ref{eq:ProjecAppQuad}).

\subsection{Numerical verification}

We proceed to verify that the proposed approach for both linear and
quadratic elements does in fact restore optimal convergence. We do
this by solving the example problem (\ref{eq:ModelProblem1}) using
the proposed discrete bilinear form instead of (\ref{eq:a_hQ}). 
\begin{figure}
\begin{centering}
\includegraphics[scale=0.5]{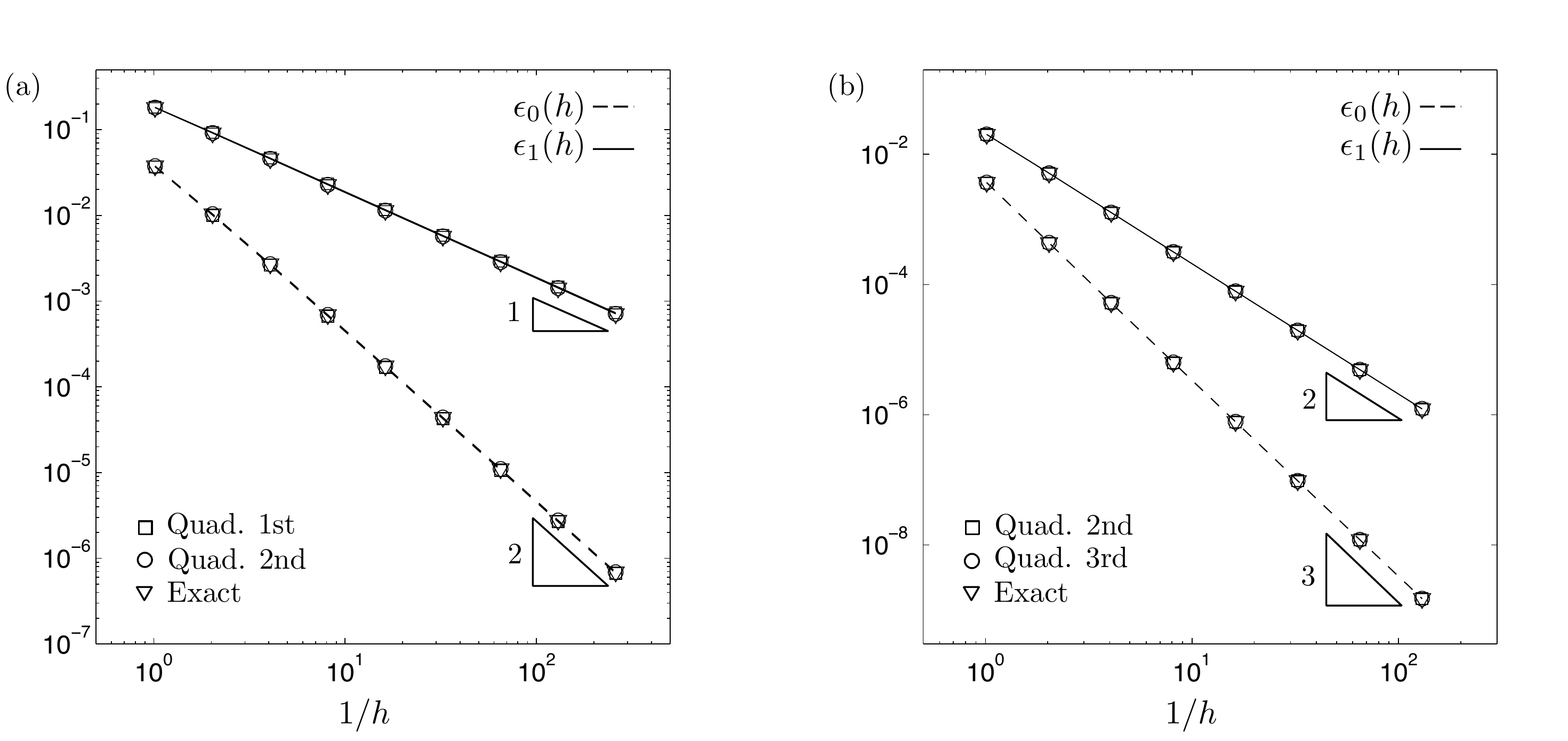}
\par\end{centering}

\caption{Results of the convergence study with proposed splitting approach
with (a) linear elements (b) quadratic elements. In both cases, the
quadrangulation scheme with indicated integration order is used \label{fig:fig3}}
\end{figure}

As shown in Figure \ref{fig:fig3}(a), we recover optimal convergence
rates for linear elements even with the first order quadrature rule.
Moreover, the solution errors are nearly identical to those obtained
from exact integration, even with the low-order quadrature. In fact,
the largest difference in the energy norm errors between the proposed
scheme with the first order rule and exact integration is $2.7\%$.
The $L^{2}$-error is slightly smaller with the proposed approach
with a difference of $4.0\%$. This shows that the first order rule,
with $n$ integration points for an $n$-gon, can be used in practice
without sacrificing accuracy. The use of the more accurate second
order error lowers this difference (to $0.027\%$ and $0.43\%$ for
the energy and $L^{2}$-norm errors, respectively) but requires four
times as many integration points.

Figure \ref{fig:fig3}(b) summarize the results for the quadratic
elements. The same conclusions can be drawn in this case: optimal
convergences rate are restored and the solution errors are almost
identical to those with the Galerkin approximation with the exact
bilinear. The largest numerical difference between the energy and
$L^{2}$-norm errors, with the second order rule, are $0.08\%$ and
$0.80\%$, respectively. 

We remark that the overhead associated with the splitting of the bilinear
form is small and the overall cost of construction of $a_{h}^{E}$
is still determined by the number of integration points. For example,
observe that (\ref{eq:LinearProjFinal}) requires visiting each edge
once, the cost of which is small compared to the geometric construction
of gradients of the Wachspress basis functions at each integration
point. In our implementation, the overhead associated with splitting
of the bilinear form (including the calculation of the projection
map) accounts for about $10\%$ of the total cost of computing the
stiffness matrix for both linear and quadratic elements. The major
difference between (\ref{eq:a_hQ}) and (\ref{eq:a_hNEW}) is that
the effort associated with basis function construction and integration
is used only on the non-polynomial part of $\mathcal{V}_{m}(E)$ where
it is needed.

\section{Treatment of non-constant coefficients }

We now briefly discuss a possible strategy to handle diffusion tensors
with variable coefficients. Such treatment of a position-dependent
material function is relevant for a wide range of problems such heat
transfer in systems with variable thermal conductivity, flow in porous
media with variable permeability, electric conduction with variable
resistivity, and magnetostatics with variable magnetic permeability
\citep{Sutradhar:2004p3652,Paulino2002,Paulino2004}. In this general
case, we are given a symmetric tensor $\mathbb{K}\in L^{\infty}(\Omega)^{2\times2}$
satisfying
\begin{equation}
\alpha^{-1}\left|\boldsymbol{\xi}\right|^{2}\leq\boldsymbol{\xi}\cdot\mathbb{K}(\boldsymbol{x})\boldsymbol{\xi}\leq\alpha\left|\boldsymbol{\xi}\right|^{2},\qquad\forall\boldsymbol{\xi}\in\mathbb{R}^{2},\ \forall\boldsymbol{x}\in\Omega\label{eq:K-coerc}
\end{equation}
for some positive constant $\alpha$. In order for use of quadrature
to make sense, we assume that $\mathbb{K}$ is defined everywhere
in the domain.

For each element $E\in\mathcal{T}_{h}$, we construct a first-order
approximation to $\mathbb{K}$ over $E$, which we denote by $\mathbb{K}_{E}$.
For example, we can do so by setting
\begin{equation}
\mathbb{K}_{E}=\frac{1}{\left|E\right|}\int_{E}\mathbb{K}\mathrm{d}\boldsymbol{x}
\end{equation}
or, if $\mathbb{K}$ is a smooth function, we can take $\mathbb{K}_{E}$
to be simply the value of $\mathbb{K}$ the center of $E$. For linear
elements, we can proceed as before but with $\mathbb{K}_{E}$ in place
of $\mathbb{K}$ without sacrificing first-order convergence rate.
For quadratic elements, however, this will lead to loss of an order
of convergence. Therefore, we consider the construction of bilinear
form that includes a correction term:

\begin{eqnarray}
a_{h}^{E}(u,v) & = & \int_{E}\nabla\Pi_{m}^{E}u\cdot\mathbb{K}_{E}\nabla\Pi_{m}^{E}v\mathrm{d}\boldsymbol{x}+\fint_{E}\nabla\left(u-\Pi_{m}^{E}u\right)\cdot\mathbb{K}_{E}\nabla\left(v-\Pi{}_{m}^{E}v\right)\mathrm{d}\boldsymbol{x}\nonumber \\
 &  & \qquad\qquad+\fint_{E}\nabla u\cdot\left(\mathbb{K}-\mathbb{K}_{E}\right)\nabla v\mathrm{d}\boldsymbol{x}\label{eq:a_hCorrected}
\end{eqnarray}
Here $\Pi_{m}^{E}$ is defined as in the previous section with $\mathbb{K}$
replaced by $\mathbb{K}_{E}$. With the correction term, we capture
the variation of $\mathbb{K}$ inside the element through sampling
$\mathbb{K}-\mathbb{K}_{E}$ at the integration points. At the same
time, we retain the simplicity offered by a constant tensor in constructing
the projection map.

Analysis of the convergence of the resulting approximation can be
based on the estimate (\ref{eq:Strang}). Here we content ourselves
with numerical verification of optimal performance. Borrowing from
\citep{brezzi2005family}, we consider the problem posed on $\Omega=\left]0,1\right[^{2}$
with exact solution 
\begin{figure}
\begin{centering}
\includegraphics[scale=0.5]{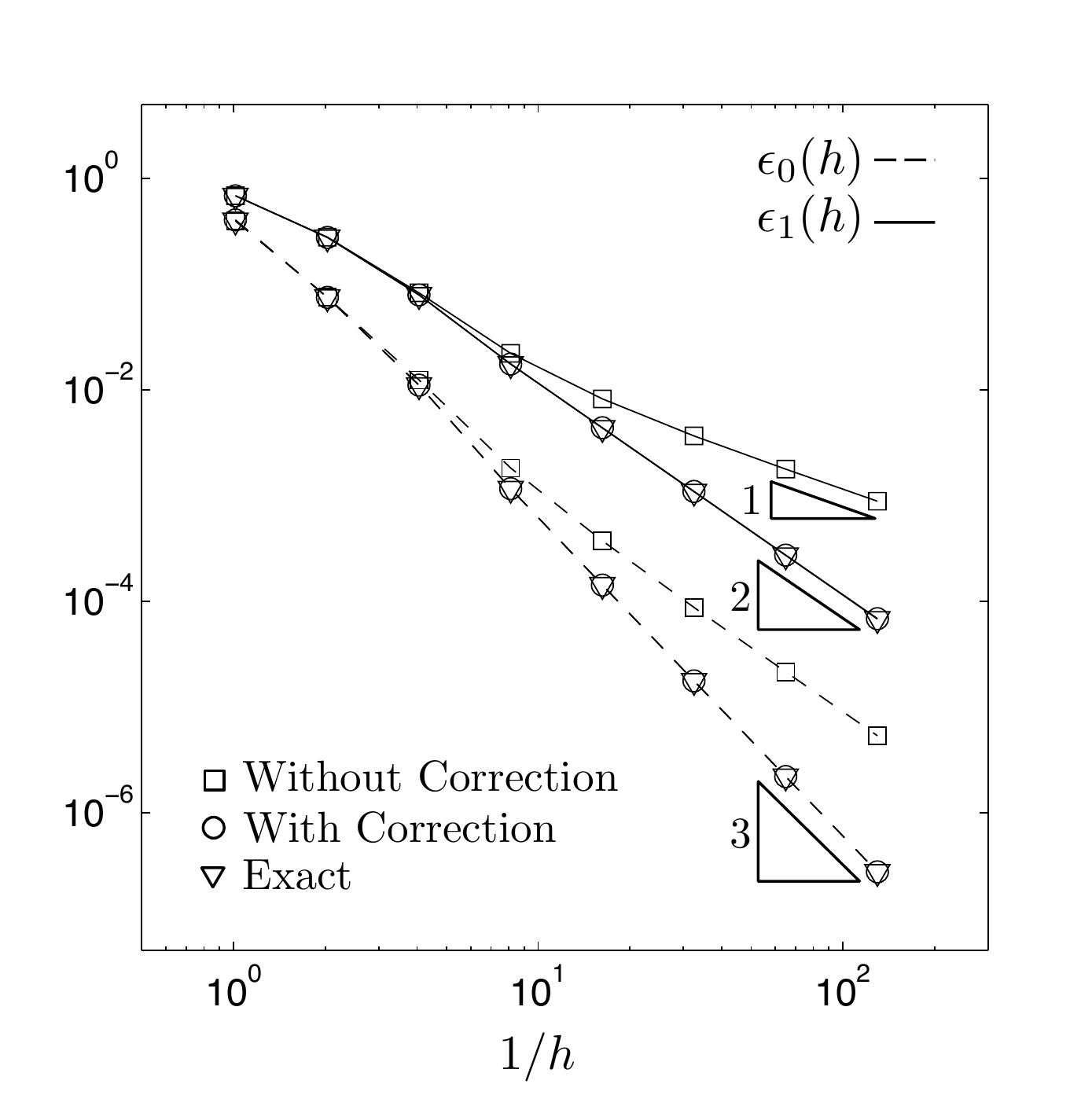}
\par\end{centering}

\caption{Results of the convergence study for problem (\ref{eq:ModelProblem2})-(\ref{eq:ModelProblem2-K})
using the quadratic elements and second-order quadrangulation scheme\label{fig:fig4}}
\end{figure}
\begin{equation}
u(\boldsymbol{x})=x_{1}^{3}x_{2}^{2}+x_{1}\sin(2\pi x_{1}x_{2})\sin(2\pi x_{2})\label{eq:ModelProblem2}
\end{equation}
and diffusion tensor given by
\begin{equation}
\mathbb{K}(\boldsymbol{x})=\left[\begin{array}{cc}
\left(x_{1}+1\right)^{2}+x_{2}^{2} & -x_{1}x_{2}\\
-x_{1}x_{2} & \left(x_{1}+1\right)^{2}
\end{array}\right]\label{eq:ModelProblem2-K}
\end{equation}
The source function $f$ and boundary data $g$ are prescribed in
accordance with (\ref{eq:ModelProblem2}) and (\ref{eq:ModelProblem2-K}).
We use the same sequence of regular meshes as before and set $\mathbb{K}_{E}=\fint_{E}\mathbb{K}\mathrm{d}\boldsymbol{x}/\left|E\right|$.
The discrete linear form uses the same quadrature that is used for
the bilinear form. We omit the results for linear elements because
the solution errors, even without the correction term, are very close
to the errors obtained with exact integration. Figure (\ref{fig:fig4})
summarizes the results for the quadratic elements. It is evident that
without the correction term, the rates of convergence are reduced
by exactly one order as a result of first-order approximation of $\mathbb{K}$.
However, the choice of (\ref{eq:a_hCorrected}) not only recovers
optimal convergence rates, but also leads to nearly the same solution
errors as for the exact integration (difference of $<1\%$).

\section{Conclusions}

We conclude by noting that the issue of quadrature error and its adverse
effect on convergence is in fact more pronounced in three dimensions.
This is because the construction of basis functions for general polyhedral
elements is more costly, quadrature rules are more difficult to obtain
and the patch test errors are typically larger. The present approach
and use of polynomial projections can help overcome the challenges
associated with polyhedral finite elements. The proposed approach
can also be beneficial for reducing the burden of integration for
nonlinear problems where integration of constitutive relations is
usually performed at the quadrature points. Finally, the formalism
of the polynomial projections is promising in furnishing an alternative
way to address the challenges of numerical integration in meshless
methods.

\section*{Acknowledgements}

This work was inspired by the ``Workshop on Discretization Methods
for Polygonal and Polyhedral Meshes'' which was held in Milan, Italy,
from September 17 to 19, 2012. In fact, the inception of this work
has its roots in our insightful discussions with Franco Brezzi. We
also appreciate all the interactions we had with the workshop participants.
We acknowledge support from the US National Science Foundation under
grant CMMI \#1321661 and from the Donald B. and Elizabeth M. Willett
endowment at the University of Illinois at Urbana-Champaign. Any opinion,
finding, conclusions or recommendations expressed here are those of
the authors and do not necessarily reflect the views of the sponsors.

\section*{Appendix A: Construction of polygonal elements}

We discuss the construction of element spaces $\mathcal{V}_{m}(E)$
by means of generalized barycentric coordinates associated with polygon
$E$. We will describe the Wachspress coordinates but note that any
other set of barycentric coordinates (e.g. Mean Value, harmonic, Sibson,
etc.) can be also used. A possibly more economical alternative in
the two-dimensional setting makes use of the usual isoparameteric
mapping \citep{Sukumar:2004p401}. As discussed in \citep{FluidsIJNMF},
the isoparametric construction in fact defines a new set of barycentric
coordinates for polygons. We emphasize that the results of this paper
applies to all resulting elements regardless of the choice of the
barycentric coordinates.

Suppose $E$ is a strictly convex $n$-gon with vertices located at
$\boldsymbol{x}_{1},\dots,\boldsymbol{x}_{n}$ oriented counter-clockwise.
The Wachspress coordinate associated with the $i$th vertex is defined
in the interior of $E$ by
\begin{equation}
\varphi_{i}(\boldsymbol{x})=\frac{w_{i}(\boldsymbol{x})}{\sum_{j=1}^{n}w_{j}(\boldsymbol{x})}
\end{equation}
with weight functions given by
\begin{equation}
w_{i}(\boldsymbol{x})=\frac{A(\boldsymbol{x}_{i-1},\boldsymbol{x}_{i},\boldsymbol{x}_{i+1})}{A(\boldsymbol{x}_{i-1},\boldsymbol{x}_{i},\boldsymbol{x})A(\boldsymbol{x}_{i},\boldsymbol{x}_{i+1},\boldsymbol{x})}
\end{equation}
Here $A\left(\boldsymbol{a},\boldsymbol{b},\boldsymbol{c}\right)$
denotes the area of the triangle with vertices located at points $\boldsymbol{a}$,
$\boldsymbol{b}$ and $\boldsymbol{c}$. We are using the convention
that $\boldsymbol{x}_{n+1}=\boldsymbol{x}_{1}$ and $\boldsymbol{x}_{0}=\boldsymbol{x}_{n}$.
It is evident that $\varphi_{i}$'s are positive functions that form
a partition of unity in $E^{\circ}$. Moreover, one can show
\begin{equation}
\sum_{i=1}^{n}\boldsymbol{x}_{i}\varphi_{i}(\boldsymbol{x})=\boldsymbol{x},\qquad\forall\boldsymbol{x}\in E^{\circ}\label{eq:Wach-x-completeness}
\end{equation}
From these properties, it follows that Wachspress functions can be
extended continuously to $\partial E$ such that \citep{Floater:2006p406}
\begin{equation}
\varphi_{i}(\boldsymbol{x})=1-\frac{\left|\boldsymbol{x}-\boldsymbol{x}_{i}\right|}{\left|\boldsymbol{x}_{i+1}-\boldsymbol{x}_{i}\right|},\qquad\varphi_{i+1}(\boldsymbol{x})=\frac{\left|\boldsymbol{x}-\boldsymbol{x}_{i}\right|}{\left|\boldsymbol{x}_{i+1}-\boldsymbol{x}_{i}\right|},\qquad\varphi_{j}(\boldsymbol{x})=0,\ \forall j\neq i,i+1\label{eq:Wach-Edge}
\end{equation}
if $\boldsymbol{x}$ lies on the edge connecting $\boldsymbol{x}_{i}$
and $\boldsymbol{x}_{i+1}$. Note that (\ref{eq:Wach-Edge}) implies
that Wachspress coordinates satisfy the Kronecker-delta property,
i.e., $\varphi_{i}(\boldsymbol{x}_{j})=\delta_{ij}$ and vary linearly
on $\partial E$%
\footnote{These properties are essential ingredients in constructing the conforming
finite element space (\ref{eq:V_h}) with degrees of freedom associated
with the vertices of the mesh. %
}. Subsequently, we have the linear precision property of
\begin{equation}
p(\boldsymbol{x})=\sum_{i=1}^{n}p(\boldsymbol{x}_{i})\varphi_{i}(\boldsymbol{x}),\qquad\forall p\in\mathcal{P}_{1}(E)\label{eq:LinearPrecision}
\end{equation}
for any point $\boldsymbol{x}$ in the closure of $E$. We set the
linear element space for $E$ as $\mathcal{V}_{1}(E)=\mbox{span}\left\{ \varphi_{1},\dots,\varphi_{n}\right\} $.
Observe how (\ref{eq:LinearPrecision}) implies (\ref{eq:VEsup PE})
with $m=1$.

To construct the quadartic serendipity element on $E$, we first define
mid-side nodes $\hat{\boldsymbol{x}}_{i}=\left(\boldsymbol{x}_{i}+\boldsymbol{x}_{i+1}\right)/2$.
The basis functions for $\mathcal{V}_{2}(E)$ are given by
\begin{equation}
\psi_{i}(\boldsymbol{x})=\sum_{a=1}^{n}\sum_{b=1}^{n}c_{i}^{ab}\varphi_{a}(\boldsymbol{x})\varphi_{b}(\boldsymbol{x}),\quad i=1,\dots,2n
\end{equation}
where $\varphi_{a}$ are barycentric coordinates for $E$ and coefficients
$c_{i}^{ab}$ are chosen such that
\begin{equation}
p(\boldsymbol{x})=\sum_{i=1}^{n}\left[p(\boldsymbol{x}_{i})\psi_{i}(\boldsymbol{x})+p(\hat{\boldsymbol{x}}_{i})\psi_{i+n}(\boldsymbol{x})\right],\qquad\forall p\in\mathcal{P}_{2}(E)
\end{equation}
and Kronecker-delta property is satisfied%
\footnote{That is, $\psi_{i}(\boldsymbol{x}_{j})=\psi_{i+n}(\hat{\boldsymbol{x}}_{j})=\delta_{ij}$
and $\psi_{i}(\hat{\boldsymbol{x}}_{j})=\psi_{i+n}(\boldsymbol{x}_{j})=0$%
}. In \citep{rand2011quadratic}, it is shown that a stable choice
of coefficients $c_{i}^{ab}$ exists and a procedure for computing
them is presented. As a result of this construction, the basis functions
exhibit quadratic variation on the boundary and (\ref{eq:VEsup PE})
is satisfied for $m=2$.

\section*{Appendix B: Implementation aspects}

We will provide details on the algebraic construction of the projection
map and the discrete bilinear forms for both linear and quadratic
elements. The presentation proceeds along similar lines as \citep{RussoSlides}
where implementation of a first order VEM formulation for Poisson's
problem is discussed.

We unify the presentation by noting that the right hand side of (\ref{eq:ProjecAppQuad})
reduces to $\int_{\partial E}v\mathbb{K}\nabla p\cdot\boldsymbol{n}\mathrm{d}s=a^{E}(v,p)$
whenever $p\in\mathcal{P}_{1}(E)$. Therefore, the first condition
in the definition of the projection maps is taken to be%
\footnote{For the linear element, this term vanishes. For the quadratic element,
as in section 5.2, we assume that the quadrature used for the first
term is exact when $v\in\mathcal{P}_{2}(E)$. %
}
\begin{equation}
a^{E}(p,\Pi_{m}^{E}v)=-\fint_{E}v\Div\left(\mathbb{K}\nabla p\right)\mathrm{d}\boldsymbol{x}+\int_{\partial E}v\mathbb{K}\nabla p\cdot\boldsymbol{n}\mathrm{d}s,\qquad\forall p\in\mathcal{P}_{m}(E)\label{eq:ProjectionApp}
\end{equation}
for both linear and quadratic elements. We can also replace the second
condition in (\ref{eq:ProjectionExact}) with an equivalent condition
given by: 
\begin{equation}
\overline{\Pi_{m}^{E}v}=\overline{v},\qquad\forall v\in\mathcal{V}_{m}(E)\label{eq:MeanValueEquality}
\end{equation}
The equivalence follows from the fact that (\ref{eq:ProjectionApp})
implies that $\nabla\Pi_{m}^{E}p=\nabla p$ for $p\in\mathcal{P}_{m}(E)$.
Thus, (\ref{eq:ProjectionApp}) together with (\ref{eq:MeanValueEquality})
ensures that the projection fixes polynomial functions.

Let $n_{v}=\dim\mathcal{V}_{m}(E)$ and $n_{p}=\dim\mathcal{P}_{m}(E)-1$,
and consider a basis for $\mathcal{P}_{m}(E)$ of the form $\left\{ 1,p_{1},\dots,p_{n_{p}}\right\} $
such that $\overline{p_{\alpha}}=0$ for $\alpha=1,\dots,n_{p}$.
For example,
\begin{equation}
p_{0}(\boldsymbol{x})=1,\quad p_{1}(\boldsymbol{x})=x_{1}-\overline{x_{1}},\quad p_{2}(\boldsymbol{x})=x_{2}-\overline{x_{2}}
\end{equation}
is such a basis for $\mathcal{P}_{1}(E)$. As before, let $\left\{ \varphi_{1},\dots,\varphi_{n_{v}}\right\} $
be the canonical basis for $\mathcal{V}_{m}(E)$. 

We define two matrices $\boldsymbol{R}$ and $\boldsymbol{N}$ of
size $n_{v}\times n_{p}$ by
\begin{eqnarray}
\boldsymbol{R}_{i\alpha} & = & -\fint_{E}\varphi_{i}\Div\left(\mathbb{K}\nabla p_{\alpha}\right)\mathrm{d}\boldsymbol{x}+\int_{\partial E}\varphi_{i}\mathbb{K}\nabla p_{\alpha}\cdot\boldsymbol{n}\mathrm{d}s\label{eq:R_matrix}\\
\boldsymbol{N}_{i\alpha} & = & p_{\alpha}(\boldsymbol{x}_{i})\label{eq:N_matrix}
\end{eqnarray}
where $\boldsymbol{x}_{i}$ is the location of the $i$th node of
$E$ (associated with $\varphi_{i}$). Observe that $\boldsymbol{R}_{i\alpha}$
is the right-hand-side of (\ref{eq:ProjectionApp}) for $v=\varphi_{i}$
and $p=p_{\alpha}$. Also the Lagrangian property of the basis functions
and their polynomial precision implies
\begin{equation}
p_{\alpha}(\boldsymbol{x})=\sum_{i=1}^{n_{v}}\boldsymbol{N}_{i\alpha}\varphi_{i}(\boldsymbol{x})\label{eq:p_alpha_expansion}
\end{equation}
Using (\ref{eq:R_matrix})-(\ref{eq:p_alpha_expansion}) and the exactness
of the quadrature rule on polynomials, we have
\begin{eqnarray}
a^{E}(p_{\alpha},p_{\beta}) & = & -\int_{E}p_{\alpha}\Div\left(\mathbb{K}\nabla p_{\beta}\right)\mathrm{d}\boldsymbol{x}+\int_{\partial E}p_{\alpha}\mathbb{K}\nabla p_{\beta}\cdot\boldsymbol{n}\mathrm{d}s\nonumber \\
 & = & -\fint_{E}p_{\alpha}\Div\left(\mathbb{K}\nabla p_{\beta}\right)\mathrm{d}\boldsymbol{x}+\int_{\partial E}p_{\alpha}\mathbb{K}\nabla p_{\beta}\cdot\boldsymbol{n}\mathrm{d}s\label{eq:NR}\\
 & = & \sum_{i=1}^{n_{v}}\boldsymbol{N}_{i\alpha}\left[-\fint_{E}\varphi_{i}\Div\left(\mathbb{K}\nabla p_{\beta}\right)\mathrm{d}\boldsymbol{x}+\int_{\partial E}\varphi_{i}\mathbb{K}\nabla p_{\beta}\cdot\boldsymbol{n}\mathrm{d}s\right]\nonumber \\
 & = & \left[\boldsymbol{N}^{T}\boldsymbol{R}\right]_{\alpha\beta}\nonumber 
\end{eqnarray}
This shows that $\boldsymbol{N}^{T}\boldsymbol{R}$ is an $n_{p}\times n_{p}$
symmetric positive definite matrix.

Since $\Pi_{m}^{E}\varphi_{i}$ is an element of $\mathcal{P}_{m}(E)$,
there exists an $n_{v}\times n_{p}$ matrix $\boldsymbol{S}$ such
that: 
\begin{equation}
\Pi_{m}^{E}\varphi_{i}=\frac{1}{n_{v}}+\sum_{\beta=1}^{n_{p}}\boldsymbol{S}_{i\beta}p_{\beta}\label{eq:Pi_phi_expansion}
\end{equation}
for $i=1,\dots,n_{v}$. Note that $\overline{\Pi_{m}^{E}\varphi_{i}}=1/n_{v}=\overline{\varphi_{i}}$
and so (\ref{eq:MeanValueEquality}) is satisfied. To derive an expression
for $\boldsymbol{S}$, we appeal to (\ref{eq:ProjectionApp}), and
set $p=p_{\alpha}$ and $v=\varphi_{i}$ to get
\begin{equation}
\sum_{\beta=1}^{n_{p}}\boldsymbol{S}_{i\beta}a^{E}(p_{\alpha},p_{\beta})=\boldsymbol{R}_{i\alpha}\label{eq:S-R_intermediate}
\end{equation}
Here we have used the linearity of the bilinear form and expansion
(\ref{eq:Pi_phi_expansion}).

From (\ref{eq:NR}) and the fact that (\ref{eq:S-R_intermediate})
must hold for $\alpha=1,\dots,n_{p}$, we obtain following expression
for $\boldsymbol{S}$ 
\begin{equation}
\boldsymbol{S}=\boldsymbol{R}\left(\boldsymbol{N}^{T}\boldsymbol{R}\right)^{-1}\label{eq:S}
\end{equation}
which in turn, through (\ref{eq:Pi_phi_expansion}), gives the expression
for the projection map. One can verify that for $m=1$, (\ref{eq:Pi_phi_expansion})
and (\ref{eq:S}) recover the expression (\ref{eq:LinearProjFinal})
derived in section 5.1.

We proceed to derive explicit algebraic expressions for the stiffness
matrix associated with bilinear $a_{h}^{E}$. Recall that the $\left(i,j\right)$th
entry of the stiffness matrix associated with $E$ is given by 
\begin{equation}
a_{h}^{E}(\varphi_{i},\varphi_{j})=a^{E}(\Pi_{m}^{E}\varphi_{i},\Pi_{m}^{E}\varphi_{j})+\fint_{E}\nabla\left(\varphi_{i}-\Pi_{m}^{E}\varphi_{i}\right)\cdot\mathbb{K}\nabla\left(\varphi_{j}-\Pi_{m}^{E}\varphi_{j}\right)\mathrm{d}\boldsymbol{x}\label{eq:StiffnessMatrixIJ}
\end{equation}
We can compute the first term of the stiffness matrix as follows
\begin{equation}
a^{E}(\Pi_{m}^{E}\varphi_{i},\Pi_{m}^{E}\varphi_{j})=\sum_{\alpha=1}^{n_{p}}\sum_{\beta=1}^{n_{p}}\boldsymbol{S}_{i\alpha}\boldsymbol{S}_{j\beta}a^{E}(p_{\alpha},p_{\beta})=\left[\boldsymbol{S}\left(\boldsymbol{N}^{T}\boldsymbol{R}\right)\boldsymbol{S}^{T}\right]_{ij}=\left[\boldsymbol{R}\left(\boldsymbol{N}^{T}\boldsymbol{R}\right)^{-1}\boldsymbol{R}^{T}\right]_{ij}
\end{equation}
To get an expression for the second term, we note that
\begin{eqnarray}
\varphi_{i}-\Pi_{m}^{E}\varphi_{i} & = & \varphi_{i}-\frac{1}{n_{v}}-\sum_{\beta=1}^{n_{p}}\boldsymbol{S}_{i\beta}\sum_{j=1}^{n_{v}}\boldsymbol{N}_{j\beta}\varphi_{j}\nonumber \\
 & = & \varphi_{i}-\frac{1}{n_{v}}-\sum_{j=1}^{n_{v}}\left(\boldsymbol{S}\boldsymbol{N}^{T}\right)_{ij}\varphi_{j}\\
 & = & \sum_{j=1}^{n_{v}}\left(\boldsymbol{I}_{ij}-\frac{1}{n_{v}}\boldsymbol{U}_{ij}-\left[\boldsymbol{R}\left(\boldsymbol{N}^{T}\boldsymbol{R}\right)^{-1}\boldsymbol{N}^{T}\right]_{ij}\right)\varphi_{j}\nonumber 
\end{eqnarray}
where $\boldsymbol{I}$ denotes the $n_{v}\times n_{v}$ identity
matrix and $\boldsymbol{U}$ is the $n_{v}\times n_{v}$ matrix with
unit entries. Defining, 
\begin{equation}
\boldsymbol{P}:=\boldsymbol{I}-\frac{1}{n_{v}}\boldsymbol{U}-\boldsymbol{R}\left(\boldsymbol{N}^{T}\boldsymbol{R}\right)^{-1}\boldsymbol{N}^{T}\label{eq:P}
\end{equation}
we have $\varphi_{i}-\Pi_{m}^{E}\varphi_{i}=\sum_{j=1}^{n_{v}}\boldsymbol{P}_{ij}\varphi_{j}$.
In turn, the second term of (\ref{eq:StiffnessMatrixIJ}) can be written
as 
\begin{equation}
\fint_{E}\nabla\left(\varphi_{i}-\Pi_{m}^{E}\varphi_{i}\right)\cdot\mathbb{K}\nabla\left(\varphi_{j}-\Pi_{m}^{E}\varphi_{j}\right)\mathrm{d}\boldsymbol{x}=\sum_{k=1}^{n_{v}}\sum_{\ell=1}^{n_{v}}\boldsymbol{P}_{ik}\boldsymbol{P}_{j\ell}\left(\fint_{E}\nabla\varphi_{k}\cdot\mathbb{K}\nabla\varphi_{\ell}\mathrm{d}\boldsymbol{x}\right)
\end{equation}
Observe that the term in the parenthesis is $(k,\ell)$th entry of
the usual stiffness matrix obtained from quadrature (i.e., the stiffness
matrix corresponding to bilinear form (\ref{eq:a_hQ})). Defining
$\boldsymbol{K}_{k\ell}:=\fint_{E}\nabla\varphi_{k}\cdot\mathbb{K}\nabla\varphi_{\ell}\mathrm{d}\boldsymbol{x}$,
the expression for the stiffness matrix reduces to 
\begin{equation}
a_{h}^{E}(\varphi_{i},\varphi_{j})=\left[\boldsymbol{R}\left(\boldsymbol{N}^{T}\boldsymbol{R}\right)^{-1}\boldsymbol{R}^{T}+\boldsymbol{P}\boldsymbol{K}\boldsymbol{P}^{T}\right]_{ij}
\end{equation}

For case of non-constant coefficients, the matrix $\boldsymbol{R}$
is defined as (\ref{eq:R_matrix}) but with $\mathbb{K}$ replaced
by $\mathbb{K}_{E}$. Setting $\tilde{\boldsymbol{K}}_{ij}=\fint_{E}\nabla\varphi_{i}\cdot\mathbb{K}\nabla\varphi_{j}\mathrm{d}\boldsymbol{x}$
and $\boldsymbol{K}_{ij}=\fint_{E}\nabla\varphi_{i}\cdot\mathbb{K}_{E}\nabla\varphi_{j}\mathrm{d}\boldsymbol{x}$,
the stiffness matrix associated with the corrected bilinear (\ref{eq:a_hCorrected})
is 
\begin{equation}
a_{h}^{E}(\varphi_{i},\varphi_{j})=\left[\boldsymbol{R}\left(\boldsymbol{N}^{T}\boldsymbol{R}\right)^{-1}\boldsymbol{R}^{T}+\boldsymbol{P}\boldsymbol{K}\boldsymbol{P}^{T}+(\tilde{\boldsymbol{K}}-\boldsymbol{K})\right]_{ij}
\end{equation}
where $\boldsymbol{P}$ is again defined by (\ref{eq:P}).

\section*{References}

\bibliographystyle{siam}
\bibliography{TalischiPaulino-IntegErrPolyFE}

\end{document}